\documentclass{amsart}
\usepackage{amsmath, amsthm, amssymb, amscd, mathrsfs, eucal, epsfig}
\usepackage{pinlabel}

\begin{document}

\newtheorem{theorem}{Theorem}[section]
\newtheorem{lemma}[theorem]{Lemma}
\newtheorem{corollary}[theorem]{Corollary}
\newtheorem{conjecture}[theorem]{Conjecture}
\newtheorem{question}[theorem]{Question}
\newtheorem{problem}[theorem]{Problem}
\newtheorem*{claim}{Claim}
\newtheorem*{criterion}{Criterion}
\newtheorem*{surgery_thm}{Theorem A}
\newtheorem*{surgery_thm2}{Theorem B}
\newtheorem*{length_gap_thm}{Theorem C}
\newtheorem*{first_limit_thm}{Theorem D}
\newtheorem*{bad_comparison_thm}{Theorem E}

\theoremstyle{definition}
\newtheorem{definition}[theorem]{Definition}
\newtheorem{construction}[theorem]{Construction}
\newtheorem{notation}[theorem]{Notation}
\newtheorem{convention}[theorem]{Convention}
\newtheorem*{warning}{Warning}

\theoremstyle{remark}
\newtheorem{remark}[theorem]{Remark}
\newtheorem{example}[theorem]{Example}

\def\cl{\textnormal{cl}}
\def\scl{\textnormal{scl}}
\def\drl{\textnormal{drl}}
\def\area{\text{area}}
\def\vol{\text{volume}}
\def\rot{\textnormal{rot}}
\def\id{\text{id}}
\def\H{\mathbb H}
\def\Z{\mathbb Z}
\def\R{\mathbb R}
\def\Q{\mathbb Q}
\def\F{\mathcal F}
\def\SL{\textnormal{SL}}
\def\PSL{\textnormal{PSL}}
\def\RP{\mathbb RP}
\def\CAT{\textnormal{CAT}}
\def\Aut{\textnormal{Aut}}
\def\Out{\textnormal{Out}}
\def\Sol{\textnormal{Sol}}
\def\Nil{\textnormal{Nil}}
\def\til{\widetilde}
\def\length{\textnormal{length}}
\def\complexlength{\textnormal{complex length}}
\def\genus{\textnormal{genus}}
\def\axis{\textnormal{axis}}

\title{Length and stable length}
\author{Danny Calegari}
\address{Department of Mathematics \\ Caltech \\ Pasadena CA, 91125}
\email{dannyc@its.caltech.edu}
\date{6/3/2007, Version 0.22}

\begin{abstract}
This paper establishes the existence of a {\em gap} for the stable
length spectrum on a hyperbolic manifold. If $M$ is a hyperbolic $n$-manifold,
for every positive $\epsilon$ there is a positive $\delta$ depending only
on $n$ and on $\epsilon$ such that an element of $\pi_1(M)$ with stable
commutator length less than $\delta$ is represented by a geodesic with length
less than $\epsilon$. Moreover, for any such $M$, the first accumulation point
for stable commutator length on conjugacy classes is at least $1/12$.

Conversely, ``most'' short geodesics in hyperbolic $3$-manifolds have arbitrarily
small stable commutator length. Thus stable commutator
length is typically good at detecting the thick-thin decomposition of $M$,
and $1/12$ can be thought of as a kind of {\em homological Margulis constant}.
\end{abstract}

\maketitle

\section{Introduction}

If $M$ is a closed hyperbolic $n$-manifold with $n>2$, the Mostow Rigidity Theorem
says that the geometric
structure on $M$ is uniquely determined by the algebraic structure of $\pi_1(M)$,
up to isometry. This means that in principle, any kind of
geometric information about $M$ (e.g. length spectrum, volume, systoles, etc.)
can be recovered from group theory. On the other hand, with some notable exceptions, it is
not easy to find concepts which have a straightforward interpretation on both
the geometric and the algebraic side.

\vskip 12pt

One of the most important geometric structures on a hyperbolic manifold is the
{\em thick-thin} decomposition. There is a universal positive constant $\epsilon(n)$ in
each dimension (i.e. the {\em Margulis constant}; see \cite{Kazhdan_Margulis}) 
such that the part of a hyperbolic
$n$-manifold $M$ with injectivity radius less than $\epsilon$ (i.e. the ``thin" piece)
consists of cusps and solid tubes. In the complement of these Margulis pieces (i.e.
the ``thick" piece) the geometry and topology is uniformly controlled by the volume. 
The existence of the Margulis constant implies that in each dimension, there is a
{\em universal} notion of what it means for a closed geodesic to be {\em short}.
A careful analysis of the behavior of thin pieces under various operations 
is the source of many compactness and finiteness theorems.

In this paper, we establish a relationship between the {\em hyperbolic length} of a
closed geodesic $\gamma$, and the {\em stable commutator length} (to be defined in the
sequel) of the corresponding conjugacy class in $\pi_1(M)$.
Parallel to the geometric relationship between hyperbolic length and Margulis' constant,
we show that there is a {\em universal} notion of what it means for a conjugacy
class to have small stable commutator length.

The results in this paper show that when the stable commutator length of
a conjugacy class in $\pi_1(M)$ is sufficiently small, 
so is the length of the closed geodesic in the corresponding free homotopy
class of loops in $M$, and (in the typical case),
conversely. Thus the thick-thin decomposition of $M$ is typically detected
by stable commutator length.

A precise definition of stable commutator length in groups is given in 
\S~\ref{scl_definition_subsection}.
Informally, if $G$ is a group and $a$ is an element of the commutator subgroup
(denoted $[G,G]$), the {\em commutator length} of $a$, denoted $\cl(a)$, is the
minimum number of commutators in $G$ whose product is equal to $a$. The
{\em stable commutator length} of $a$, denoted $\scl(a)$, is the liminf of
the ratio of $\cl(a^n)$ to $n$ as $n \to +\infty$. By convention define
$\cl(a)=\infty$ if $a$ is not in $[G,G]$, so that $\scl(a)=\infty$ if and only
if the image of $a$ in $H_1(G;\R)$ is nontrivial. Note that
$\scl(\cdot)$ is a {\em characteristic} function; i.e. it is constant on orbits
of $\Aut(G)$. In particular, it is constant on conjugacy classes in $G$.

\vskip 12pt

A detailed summary of the contents of this paper follows.

\vskip 12pt

Sections \S~\ref{scl_section} and \S~\ref{quasimorphism_section} are introductory,
and contain a brief exposition of the theory of stable commutator length and its relation
to $2$-dimensional bounded cohomology and quasimorphisms via Bavard's Duality Theorem. 
References for these sections are \cite{Bavard} or \cite{Calegari_scl_monograph}.

\vskip 12pt

Section \S~\ref{surgery_section} exhibits a relationship between stable
commutator length and Dehn surgery in hyperbolic $3$-manifolds. A fundamental
observation is that elements with very small stable commutator length may
be obtained in $3$-manifold groups by Dehn filling a cusped manifold. 

Let $M$ be a $3$-manifold with a single torus boundary $T=\partial M$.
Choose a meridian and a longitude $m,l$ for $T$ such that the curve $l$
is a longitude in $T$; i.e. its image in $H_1(M;\R)$ is trivial. Let
$M(p/q)$ denote the result of $p/q$ surgery on $M$. That is, $M(p/q)$ is
obtained from $M$ by Dehn filling $M$ along the curve $pm + ql$ in $T$.
If $M(p/q)$ is hyperbolic, let $\gamma(p/q)$ denote the
core geodesic of the filling torus.

Since the image of $l$ in $H_1(M;\R)$ is trivial, some multiple of $l$ bounds
a surface in $M$. Let $S$ be an orientable surface in $M$ 
which bounds $n$ times the longitude, and
let $\chi_\Q = \chi(S)/n$. If $a(p/q) \in \pi_1(M(p/q))$ is in the conjugacy class
corresponding to $\gamma(p/q)$, there is an inequality 
$$\scl(a(p/q)) \le \frac {-\chi_\Q} {2p}$$
From this the following estimates can be derived.

\begin{surgery_thm}
Let $M(p/q)$ be a hyperbolic $3$-manifold, obtained by $p/q$
surgery on a $1$-cusped hyperbolic $3$-manifold $M$. When the core geodesic
$\gamma(p/q) \subset M(p/q)$ is contained in an embedded tube of radius $T$ at
least $2$ then
$$\length(\gamma(p/q)) \le \left( \frac {3.993\pi\chi_\Q (T+1)} {Tp} \right)^2$$
\end{surgery_thm}

\begin{surgery_thm2}
Let $M$ be a $1$-cusped hyperbolic $3$-manifold. Normalize the Euclidean structure
on the cusp to have area $1$, and let $m$ be the shortest curve on
the cusp which is homologically essential in $M$. Then
$$-\chi_\Q \ge \frac 1 {2\pi\;\length(m)^2}$$
\end{surgery_thm2}

Similar estimates are implicit in the work of Agol \cite{Agol}.

From a theorem of Neumann-Zagier \cite{Neumann_Zagier}, 
it follows that these estimates are of the
correct order of magnitude.
The results in \S~\ref{surgery_section} may be summarized by saying that one can estimate the
length of ``most'' sufficiently short geodesics in a hyperbolic $3$-manifold
(up to an order of magnitude) 
and the geometry of a Margulis tube directly from stable commutator length.

\vskip 12pt

Section \S~\ref{spectral_gap_section} exhibits a relationship between
stable commutator length and geodesic length which is converse to the
relationship studied in \S~\ref{surgery_section}.

The following theorem shows that an element with sufficiently small
stable commutator length is represented by an arbitrarily short geodesic.

\begin{length_gap_thm}[Length inequality]
For every dimension $n$ and any $\epsilon > 0$ there is a constant $\delta(\epsilon,n)>0$ 
such that if $M$ is a complete hyperbolic $n$-manifold, and $a \in \pi_1(M)$
has stable commutator length $\le \delta(\epsilon,n)$ then if $\gamma$ is the
unique geodesic in the free homotopy class associated to the conjugacy class of $a$,
$$\length(\gamma) \le \epsilon$$
\end{length_gap_thm}

The dependence of $\delta$ on $\epsilon$ is not proper; as $\epsilon \to \infty$
the function $\delta(\epsilon,n)$ approaches some limit. A uniform estimate for
this upper bound, valid in any dimension, can be obtained, leading to
a {\em universal} estimate on the first accumulation point for stable commutator
length on conjugacy classes in a hyperbolic manifold.

\begin{first_limit_thm}[Spectral Gap]
Let $M$ be a closed hyperbolic manifold, of any dimension. 
Let $\delta_\infty$ be
the first accumulation point for stable commutator length on conjugacy classes. Then 
$$\frac 1 {12} \le \delta_\infty \le \frac 1 2$$
\end{first_limit_thm}

This theorem implies that for any closed hyperbolic $n$-manifold, and any
$\delta < 1/12$, there are only finitely many conjugacy classes of elements in $\pi_1$ with 
stable commutator length less than $\delta$. This gives a precise sense in
which a stable commutator length less than $1/12$ should be thought of as
{\em small}. The lower bound of $1/12$
is the homological analogue of Margulis' constant. In contrast with the
geometric Margulis constant, this homological constant can be estimated accurately,
and is independent of dimension.
The phenomena revealed in these theorems are quite robust; 
in a forthcoming paper with Koji Fujiwara \cite{Calegari_Fujiwara}
we prove similar theorems for word-hyperbolic groups and certain
groups acting on $\delta$-hyperbolic spaces, such as amalgamated free products
and mapping class groups.

\vskip 12pt

Section \S~\ref{de_Rham_section} gives an application of 
the results in \S~\ref{spectral_gap_section}. Given an element 
$a\in \pi_1(M)$ where $M$ is hyperbolic, let $\gamma$ denote the geodesic
in the free homotopy class corresponding to the conjugacy class of $a$.
Define the {\em de Rham length} $\drl(a)$ by the formula
$\drl(a) = \sup_{\alpha} (\int_\gamma \alpha)/2\pi\|d\alpha\|$
where the supremum is taken over all smooth $1$-forms $\alpha$ on $M$, and
$\|d\alpha\|$ denotes the operator norm of the $2$-form $d\alpha$. From
results in \S~\ref{quasimorphism_section} and \S~\ref{surgery_section} there is
an inequality $\scl(a) \ge \drl(a)$. Moreover, Theorem A and Theorem B
imply that the values of $\scl(a)$ and $\drl(a)$ are usually comparable when $\gamma$
is a short geodesic.

By contrast, Theorem C together with a construction due to Thurston shows
that the values of $\scl$ and $\drl$ are not always comparable when
$\length(\gamma)$ is large.

\begin{bad_comparison_thm}
Let $M$ be a finite volume complete hyperbolic $n$-manifold with $n>2$. Then for
any $\epsilon >0$ there are infinitely many conjugacy classes $a \in \pi_1(M)$ with
$$\drl(a) \le \epsilon\; \scl(a)$$
\end{bad_comparison_thm}

\vskip 12pt

Finally, in \S~\ref{geometry_section} we briefly discuss $\scl$ in $3$-manifolds
with other (non-hyperbolic) geometric structures.

\subsection{Acknowledgments}
I would like to thank Ian Agol, Iain Aitchison,
Marc Culler, Nathan Dunfield, Koji Fujiwara, Jason Manning, Shigenori
Matsumoto and the anonymous referee for some useful comments on
the material in this paper. 
I would also like to thank Bill Thurston for some
substantial comments, for a number of key ideas, 
and for drawing my attention to \cite{Thurston_III}, 
especially Theorem~\ref{injective_wrap_estimate}, 
and its relevance to stable commutator length.

While writing this paper, I was partially supported by a Sloan Research Fellowship, and
NSF grant DMS-0405491.

\section{Stable commutator length}\label{scl_section}

The material in this section and in \S~\ref{quasimorphism_section} 
is standard, and is included for the benefit
of the reader. A basic reference is \cite{Bavard} or \cite{Calegari_scl_monograph}.

\subsection{Definition}\label{scl_definition_subsection}

\begin{definition}
Let $G$ be a group, and $a \in G$ an element. Suppose $a$ is in the commutator
subgroup $a \in [G,G]$. The {\em commutator length} of
$a$, denoted $\cl(a)$, is the word length of $a$ in the (usually infinite)
generating set of $[G,G]$ consisting of all commutators.
\end{definition}

In words, $\cl(a)$ is smallest number of commutators whose product is $a$.
The function $\cl(\cdot)$ can be extended to all of $G$ by setting
$\cl(a)=\infty$ by convention when $a$ is not in $[G,G]$.

\begin{definition}
The {\em stable commutator length} of $a$, denoted $\scl(a)$ is equal to
$$\scl(a) = \liminf_{n \to \infty} \frac {\cl(a^n)} n$$
\end{definition}

The value of $\scl(a)$ is infinite if and only if the image of $a$ in $H_1(G;\Z) = G/[G,G]$
has infinite order. The image of a commutator under a homomorphism between groups
is again a commutator. It follows that the functions $\cl(\cdot)$ and $\scl(\cdot)$ are 
non-increasing under homomorphisms between groups, and 
characteristic; i.e. constant on orbits of $\Aut(G)$.

\subsection{Examples}

\begin{example}[Amenable groups]\label{amenable_example_vanishes}
If $G$ is amenable, stable commutator length is identically zero on $[G,G]$.
See Bavard \cite{Bavard} Cor. 2, p. 139.
\end{example}

\begin{example}[Free groups and surface groups]\label{free_example_bound}
If $F$ is free, $\scl(a)\ge 1/2$ for all nontrivial $a \in F$.
This follows from work of Duncan-Howie \cite{Duncan_Howie} Thm. 3.3, as explained
in \cite{Comerford_Edmunds}. The same arguments apply also to fundamental groups
of closed orientable surfaces of genus at least $2$, so $1/2$ is a lower
bound for $\scl$ on nontrivial elements in these groups too.

A similar estimate $\scl(a)\ge 1/6$ was obtained earlier by Culler \cite{Culler}.
\end{example}

\begin{example}[Positive twists in mapping class groups]
Let $S$ be a closed orientable surface of genus $g\ge 2$. If $a$ in the
mapping class group of $S$ is a product of $k$ right-handed Dehn twists along
essential disjoint simple closed curves in $S$ then
$$\scl(a) \ge \frac k {6(3g-1)}$$
See Korkmaz \cite{Korkmaz} Thm. 2.1 for the case $k=1$ or 
Kotschick \cite{Kotschick} for the general case.
\end{example}

\section{Quasimorphisms}\label{quasimorphism_section}

\subsection{Definition}

Informally, a quasimorphism is a real-valued function on a group which
is additive up to bounded error.

\begin{definition}
Let $G$ be a group. A {\em quasimorphism} on $G$ is a function
$$\phi:G \to \R$$
for which there exists some least real number $D(\phi)$
called the {\em defect} such that for all $a,b \in G$ there is an
inequality
$$|\phi(a) + \phi(b) - \phi(ab)| \le D(\phi)$$
\end{definition}

\begin{remark}
The number $D(\phi)$ is sometimes called the {\em error} or the
{\em deficit} of $\phi$ ({\em d\'efaut} in French).
\end{remark}

\begin{lemma}\label{product_inequality}
Let $S$ be a (possibly infinite) generating set for $G$. Let $w$ be
a word in the generators, representing an element of $G$. let $|w|$
denote the length of $w$ in the generator set, and let $w(i)$ denote
the $i$th letter of $w$. There is an inequality
$$|\; \phi(w) - \sum_{i=1}^{|w|} \phi(w(i))\; | \le (|w|-1)D(\phi)$$
\end{lemma}
\begin{proof}
If $w=vs$ where $s\in S$ then $|\phi(w) - \phi(v) - \phi(s)| \le D(\phi)$
by the definition of $D(\phi)$. Now induct on the length of $w$ and
apply the triangle inequality.
\end{proof}

\begin{definition}
A quasimorphism $\phi$ on $G$ is {\em symmetric} if $\phi(a^{-1}) = - \phi(a)$ for
all $a \in G$.
\end{definition}

For any quasimorphism $\phi$ the linear combination
$\phi'(a):= \frac 1 2 (\phi(a) - \phi(a^{-1}))$ is called the
{\em symmetrization} of $\phi$. 

\begin{lemma}\label{symmetrize_lemma}
For any quasimorphism $\phi$, the symmetrization $\phi'$ is a symmetric quasimorphism
and satisfies
$D(\phi') \le D(\phi)$.
\end{lemma}
\begin{proof}
For any $a,b$, there is an equality
$$\phi'(ab) - \phi'(a) - \phi'(b) = \frac 1 2 (\phi(ab) - \phi(a) - \phi(b)) -
\frac 1 2 (\phi(b^{-1}a^{-1}) - \phi(a^{-1}) - \phi(b^{-1}))$$
\end{proof}

\begin{definition}
A quasimorphism $\phi$ on $G$ is {\em homogeneous} if $\phi(a^n) = n\phi(a)$
for all $a \in G$ and all integers $n$.
\end{definition}

For any quasimorphism $\phi$ the limit
$$\overline{\phi}(a): = \lim_{n\to +\infty} \frac {\phi(a^n)} n$$
exists and is called the {\em homogenization} of $\phi$. 

\begin{lemma}\label{homogenize_lemma}
For any quasimorphism $\phi$, the homogenization $\overline{\phi}$ is
a homogeneous quasimorphism and satisfies
$D(\overline{\phi}) \le 2D(\phi)$.
\end{lemma}
See Bavard \cite{Bavard} Lem. 3.6, p. 142
or \cite{Calegari_scl_monograph} for a proof.

\begin{lemma}\label{conjugacy}
A homogeneous quasimorphism is symmetric, and
is constant on conjugacy classes.
\end{lemma}
\begin{proof}
Symmetry is part of the definition. For any $a,b \in G$, any homogeneous
quasimorphism $\phi$ 
and any integer $n$,
$$n|(\phi(b^{-1}ab) - \phi(a)| = |\phi(b^{-1}a^nb) -\phi(a^n)| \le 2D(\phi)$$
where the last inequality follows from the symmetry of $\phi$ and
Lemma~\ref{product_inequality}.
\end{proof}

\begin{notation}
The set of quasimorphisms on $G$ is denoted $\hat{Q}(G)$ and the
set of homogeneous quasimorphisms is denoted $Q(G)$.
\end{notation}

\subsection{Banach norm and exact sequence}

For any group $G$, the sets $\hat{Q}(G), Q(G)$ are vector subspaces of $\R^G$.
The function $D(\cdot)$ defines a pseudo-norm on these vector spaces, and vanishes
exactly on the set of quasimorphisms with defect $0$.

A quasimorphism with defect $0$ is a homomorphism to $\R$; the set of
such homomorphisms is naturally isomorphic to $H^1(G;\R)$.
The subspace $H^1(G;\R)$ is closed in $\hat{Q}(G)$,
the quotient space $\hat{Q}(G)/H^1(G;\R)$ is a normed
vector space with norm $D(\cdot)$, and $Q(G)/H^1(G;\R)$ is a
relatively closed normed vector subspace.

A function on $G$ such as an element $\phi \in \hat{Q}(G)$ extends by
linearity to define a $1$-cochain. The coboundary of the $1$-cochain
$\phi$ is denoted $\delta\phi$ and is a $2$-cocycle on $G$.
The defining property of a quasimorphism implies that $\delta\phi$ is
a bounded $2$-cocycle with norm equal to $D(\phi)$. This defines
a continuous map $\delta:\hat{Q}(G) \to H^2_b(G;\R)$.

\begin{theorem}
There is an exact sequence
$$0 \to H^1(G;\R) \to Q(G) \to H^2_b(G;\R) \to H^2(G;\R)$$
\end{theorem}

See Bavard \cite{Bavard} Prop. 3.3.1, p. 135 or
Besson \cite{Besson} for a proof and \cite{Gromov_volume}
for a more substantial discussion of bounded cohomology.

\subsection{Bavard's duality theorem}

There is a fundamental relationship between stable commutator length
and quasimorphisms in any group.

\begin{theorem}[Bavard's Duality Theorem \cite{Bavard}, p. 111]\label{Bavard_theorem}
Let $G$ be any group. Then for any $a \in [G,G]$ there is an equality
$$\scl(a) = \frac 1 2 \sup_{\phi \in Q(G)} \frac {|\phi(a)|} {D(\phi)}$$
\end{theorem}

The proof is non-constructive. It uses the Hahn-Banach theorem, $L^1$--$L^\infty$
duality, and the fact that $H_2(G;\R)$ is generated efficiently (in the $L^1$ sense)
by surfaces for which one knows the exact values of the $L^1$ norm on homology classes.

One direction of the theorem is elementary. A homogeneous quasimorphism
is constant on conjugacy classes by Lemma~\ref{conjugacy}. If
$a=[a_1,b_1]\cdots[a_n,b_n]$ for some $a_i,b_i \in G$ then 
$$|\phi(a)| = |\phi([a_1,b_1]\cdots[a_n,b_n]) - 
\sum_i\phi(a_ib_ia_i^{-1}) - \sum_i \phi(b_i^{-1})| \le (2n-1)D(\phi)$$
by symmetry and Lemma~\ref{product_inequality}. It is this direction of the
theorem that is most relevant to this paper.

\subsection{Examples of quasimorphisms}

\begin{example}[Poincar\'e \cite{Poincare}]\label{rotation_example}
Let $G$ be a group of orientation-preserving homeomorphisms of $S^1$ and
let $\hat{G}$ be the group of orientation-preserving homeomorphisms of $\R$
which cover elements of $G$ under the covering map $\R \to S^1$.
For $a \in \hat{G}$ define $\rot$ by
$$\rot(a):= \lim_{n \to +\infty} \frac {a^n(0)} n$$
Then $\rot$ is an element of $Q(\hat{G})$ with defect at most $1$.
\end{example}

\begin{example}[Brooks \cite{Brooks}]
Let $F$ be a free group on a generating set $S$. Let $w$ be a reduced word in the generators
$S$ of length at least $2$. For any reduced word $a \in F$, define
$$\phi_w(a) = \# \text{ of disjoint copies of } w \text{ in } a - \# \text{ of disjoint copies of } w^{-1}
\text{ in } a$$
Then one may verify directly that $\phi_w$ is a nontrivial quasimorphism with
defect at most $3$.
\end{example}

More examples are given in \S~\ref{derham_subsection}.

\section{Hyperbolic $3$-manifolds and Dehn surgery}\label{surgery_section}

This section relates stable commutator length to hyperbolic
geometry, especially in $3$-manifolds. The discussion is concrete and only
makes use of the ``easy'' direction of Theorem~\ref{Bavard_theorem}.

\subsection{Area inequality}

Let $\H^n$ denote hyperbolic $n$-space. Let $\Delta$ be a (possibly ideal)
geodesic triangle with
vertices $a,b,c$ and respective interior angles $\alpha,\beta,\gamma$ in $\H^n$.
Then by the Gauss-Bonnet formula, the area of $\Delta$ satisfies
$$\area(\Delta) = \pi - \alpha - \beta -\gamma$$
In particular, $\area(\Delta) \le \pi$ with equality if and only if $\Delta$ is
ideal.

\subsection{de Rham quasimorphisms}\label{derham_subsection}

The following construction is due to Barge-Ghys \cite{Barge_Ghys}:

Let $M$ be a hyperbolic manifold. Let $p \in M$ be a basepoint.
Let $G = \pi_1(M,p)$, and let $\alpha \in \Omega^1(M)$ be a $1$-form.
For each element $a \in G$ let $L_a$ denote the oriented
geodesic, based at $p$, in the  homotopy class of $a$.
For each $a \in G$, define $\phi(a)$ by
$$\phi(a) = \int_{L_a} \alpha$$

\begin{lemma}\label{form_is_quasimorphism}
Let $\alpha,p,\phi$ be as above.
The function $\phi$ is a symmetric quasimorphism on $G$ with 
$$D(\phi) \le \pi\|d\alpha\|$$
where $\|\cdot\|$ denotes the operator norm on the $2$-form $d\alpha$.
\end{lemma}
\begin{proof}
If $a,b \in G$ then the geodesics $L_a,L_b,L_{b^{-1}a^{-1}}$ form the three
sides of a geodesic triangle $\Delta$. Observe that $L_{b^{-1}a^{-1}}$ and
$L_{ab}$ are the same geodesic with opposite orientations. By
Stokes' theorem we get
$$\phi(a) + \phi(b) - \phi(ab) = \int_{\partial\Delta} \alpha = 
\int_\Delta d\alpha \le \|d\alpha\|\;\area(\Delta) \le \pi\|d\alpha\|$$
\end{proof}

Quasimorphisms $\phi$ of this kind are sometimes called {\em Barge-Ghys}
quasimorphisms, or {\em de Rham} quasimorphisms.

By Lemma~\ref{homogenize_lemma}, 
the homogenization $\overline{\phi}$ satisfies $D(\overline{\phi}) \le 2\pi\|d\alpha\|$.
By Theorem~\ref{Bavard_theorem}, for any $a\in G$ there is an inequality
$\scl(a) \ge \overline{\phi}(a)/4\pi\|d\alpha\|$. 

\begin{remark}
The homogenization $\overline{\phi}$ is obtained by setting
$\overline{\phi}(a) = \int_{\gamma_a} \alpha$ where $\gamma_a$ is the
{\em free} geodesic corresponding to the conjugacy class of $a$. We will not
use this fact in this paper.
\end{remark}

Given an element $a$, let $\gamma_a$ denote the free geodesic in the 
free homotopy class corresponding to the conjugacy
class of $a$. Observe if $p \in \gamma_a$ then $L_a= \gamma_a$ and
$\overline{\phi}(a) = \phi(a)$.

If $\gamma_a$ is a short simple geodesic contained in an embedded tube of large radius,
one can obtain lower bounds on $\phi(a)$.

\begin{lemma}\label{tube_estimate}
Let $a \in \pi_1(M)$ and let the conjugacy class of $a$ be
represented by a simple geodesic loop $\gamma$ in $M$ which is the core
of an embedded solid torus of radius $T$. Then there
is a homogeneous quasimorphism $\phi \in Q(\pi_1(M))$ with $D(\phi) \le 2\pi$ satisfying
$$\phi(a) = \length(\gamma) \sinh(T)\frac {T} {T+1}$$
\end{lemma}
\begin{proof}
Let $S$ be the solid torus of radius $T$ about $\gamma$. Let $r:S \to \R$ 
be the function whose value at every point is the hyperbolic distance to $\gamma$.
Denote the radial projection map by
$p:S \to \gamma$.

Parameterize $\gamma$ by $\theta$, so that
$d\theta$ is the length form on $\gamma$, and $\int_\gamma d\theta = \length(\gamma)$.
Pulling back by $p$ extends $\theta$ and $d\theta$ to all of $S$. Note that
$\theta$ takes values in $\R/\langle\length(\gamma)\rangle$ but $d\theta$ is a real-valued
$1$-form.

Define
$$\alpha = d\theta\cdot(\sinh(T) - \sinh(r))$$
on $S$, and extend it by $0$ outside $S$.  There is an equality
$$\|d\theta\| = 1/\cosh(r)$$
on $S$. By direct calculation, $d\alpha = \cosh(r)d\theta \wedge dr$ on $S$,
so that $\|d\alpha\|=1$.

The form $\alpha$ is not smooth along
$\partial S$, but it is Lipschitz. Let $\beta_\epsilon(r)$ be a $C^\infty$ function
on $[0,T]$ taking the value $1$ in a neighborhood of $0$ and the value $0$
in a neighborhood of $r$, and with $|\beta_\epsilon'|<1/(T-\epsilon)$ throughout, for
some small $\epsilon$.
The product $\alpha_\epsilon := \beta_\epsilon(r)\alpha$ is $C^\infty$ and satisfies
$$d\alpha_\epsilon = d\theta \wedge dr (\beta_\epsilon(r)\cosh(r) + \beta_\epsilon'(r)\sinh(r))$$
so $\|d\alpha_\epsilon\| \le 1 + 1/(T-\epsilon)$.

Let $\phi_\epsilon$ be a homogeneous quasimorphism, obtained from $\alpha_\epsilon$
by integrating along geodesics as in Lemma~\ref{form_is_quasimorphism}, and then
homogenizing. The forms $\alpha_\epsilon$ converge in $C^1$
as $\epsilon \to 0$ and therefore $\phi_\epsilon$ converges
in $Q(G)$ to some $\phi$. Multiplying $\phi$ by $T/(T+1)$ we obtain
a homogeneous quasimorphism with defect at most equal to $2\pi$, whose value
on $a$ is $\length(\gamma)\sinh(T)\cdot T/(T+1)$.
\end{proof}

\subsection{Estimates in embedded tubes}

To apply Lemma~\ref{tube_estimate} to a sufficiently short geodesic, 
it is necessary to estimate the radius of
an embedded tube in terms of the length of the core geodesic. The following 
estimate is due to Hodgson-Kerckhoff:

\begin{lemma}[Hodgson-Kerckhoff \cite{Hodgson_Kerckhoff} p.~403]\label{Hodgson_Kerckhoff_estimate}
Let $S$ be an embedded tube in a hyperbolic $3$-manifold. Let $T$ be the radius of
$S$ and $\length(\gamma)$ the length of the core geodesic. Then there is an estimate
$$\length(\gamma) \ge 0.5404 \frac {\tanh T} {\cosh(2T)} $$
\end{lemma}

See \cite{Hodgson_Kerckhoff} for a proof. Note for $T$ sufficiently large, 
we get $e^T \ge 1.03 \; \length(\gamma)^{-1/2}$.

\begin{remark}\label{Reznikov_remark}
For any dimension $n$ there exists
a constant $C_n$ such that $e^T \ge C_n \length(\gamma)^{-2/(n+1)}$.
See Reznikov \cite{Reznikov} p. 478 for a proof.
\end{remark}

Putting this together with Lemma~\ref{tube_estimate} and the
easy direction of Theorem~\ref{Bavard_theorem}, one obtains
a lower bound on the stable commutator length of $a$ in terms of the
geometry of a geodesic representative.

\vskip 12pt

Conversely, one can sometimes obtain an upper bound on stable
commutator length from topology. 
Let $M$ be a $3$-manifold with a single torus boundary $T=\partial M$.
Choose a meridian and a longitude $m,l$ for $T$ such that the curve $l$
is a longitude in $T$; i.e. its image in $H_1(M;\R)$ is trivial. Let
$M(p/q)$ denote the result of $p/q$ surgery on $M$. That is, $M(p/q)$ is
obtained from $M$ by Dehn filling $M$ along the curve $pm + ql$ in $T$,
where the additive notation denotes a particular primitive class
in $H_1(T;\Z)$ and therefore corresponds to a unique isotopy class of simple closed curve.
If $M(p/q)$ is hyperbolic, let $\gamma(p/q)$ denote the
core geodesic of the filling torus.

Since the image of $l$ in $H_1(M;\R)$ is trivial, there is some $n$ so that
$nl$ (with additive notation) is trivial in $H_1(M;\Z)$. This implies that there
is an orientable surface $S$ in $M$ which bounds $n$ times the longitude.
Let $S$ be such a surface of least genus, and for such
a surface, let $\chi_\Q = \chi(S)/n$. 

\begin{lemma}\label{Seifert_estimate_scl}
With notation as above, if $a(p/q) \in \pi_1(M(p/q))$ is in the conjugacy class
corresponding to $\gamma(p/q)$, there is an inequality 
$$\scl(a(p/q)) \le \frac {-\chi_\Q} {2p}$$
\end{lemma}
\begin{proof}
Let $S'$ be a finite cover of $S$ with two boundary components, for which
the degree of the covering $S' \to S$ is $2N$. Such a cover exists for every
integer $N>1$. The composition of $S' \to S$ with the map $S \to M$ wraps
each boundary component $nN$ times around $l$. Add a rectangle to $S'$
to boundary connect sum the two boundary components together, and map this
rectangle to a point in $l$. The result is a surface $S''$ with 
$\chi(S'') = \chi(S')-1$. Since Euler characteristic is multiplicative
under coverings, $\chi(S') = 2N\chi(S)$. The boundary component of $S''$
maps $2nN$ times around $l$ so 
$$\scl(l) \le \frac {\genus(S'')} {2nN} = \frac {3 - 2N\chi(S)} {4nN}$$
taking the limit $N \to +\infty$ gives the estimate $\scl(l) \le -\chi_\Q/2$.

Under $p/q$ Dehn surgery, the longitude $l$ is wrapped $p$ times around the core
geodesic of the filling torus. If $a(p/q)$ (hereafter $a$ for brevity)
is in this conjugacy class, then $a^p$ is the image of $l$. Since $\scl$
cannot increase under homomorphisms, $\scl(a) \le \scl(l)/p$.
\end{proof}

Lemma~\ref{Seifert_estimate_scl} 
says that away from finitely many lines in Dehn surgery space,
the core geodesics obtained by hyperbolic Dehn surgery correspond to
conjugacy classes in $\pi_1$ with
arbitrarily small stable commutator length. Since all sufficiently
short geodesics arise by geometric Dehn surgery on cusped manifolds, this
strongly suggests that ``typical'' short geodesics have very small stable commutator
length. 

Using Lemma~\ref{Hodgson_Kerckhoff_estimate} and Lemma~\ref{tube_estimate},
we can make the relation between $\scl$ and length quantitative. 
For simplicity, we treat the case of knot complements in $S^3$.

\begin{surgery_thm}
Let $M(p/q)$ be a hyperbolic $3$-manifold, obtained by $p/q$
surgery on a $1$-cusped hyperbolic $3$-manifold $M$. When the core geodesic
$\gamma(p/q) \subset M(p/q)$ is contained in an embedded tube of radius $T$ at
least $2$ then
$$\length(\gamma(p/q)) \le \left( \frac {3.993\pi\chi_\Q (T+1)} {Tp} \right)^2$$
\end{surgery_thm}
\begin{proof}
For the sake of legibility, abbreviate $\gamma(p/q)$ to $\gamma$ throughout
this proof.

By Lemma~\ref{Hodgson_Kerckhoff_estimate} and Lemma~\ref{tube_estimate}
there is a homogeneous quasimorphism $\phi$ with $D(\phi)\le 2\pi$
and satisfying 
$$\phi(a) = \length(\gamma)\sinh(T)\frac T {T+1}$$ where
$\length(\gamma) \ge 0.5404 \tanh(T)/\cosh(2T)$.
By Theorem~\ref{Bavard_theorem} and Lemma~\ref{Seifert_estimate_scl} there is an estimate
$$\frac {-\chi_\Q} {2p} \ge \scl(a) \ge \frac {\phi(a)} {2D(\phi)}$$ 
and therefore
$$\frac {-4\pi\chi_\Q} {2p} \ge \length(\gamma)\sinh(T) \frac T {T+1}$$

For $T>2$ there is an inequality
$$1.0376 \; e^{2T} \ge 2\cosh(2T)/\tanh(T)$$
and therefore
$$e^T \ge 1.0206 \; \length(\gamma)^{-1/2}$$
Similarly, for $T>2$ there is an inequality $2\sinh(T)>0.9816 \; e^T$ so
$$\frac {-4\pi\chi_\Q} {2p} \ge \length(\gamma)\sinh(T)\frac T {T+1} > 0.5009 \; \length(\gamma)^{1/2} \frac T {T+1}$$
\end{proof}

In \cite{Neumann_Zagier}, Neumann and Zagier define the following quadratic
form $Q$:
$$Q(p,q) = (\text{length of } pm + ql)^2/\area(\partial S)$$
where $\partial S$ is the horotorus boundary of the cusp on $M$, and 
$pm + ql$ is a straight curve on the horotorus (in the intrinsic Euclidean metric)
representing $p$ times the meridian plus $q$ times the longitude.
Equivalently, we can scale the Euclidean cusp to have area $1$, and the form
reduces to $Q(p,q) = (\text{length of }pm + ql)^2$.

Proposition 4.3 from \cite{Neumann_Zagier} gives the following asymptotic formula for
$\length(\gamma)$, the length of the core geodesic of $M(p/q)$:
\begin{lemma}[Neumann--Zagier]\label{Neumann_Zagier_formula}
$$\length(\gamma(p/q)) = 2\pi Q(p,q)^{-1} + O\left( \frac 1 {p^4 + q^4} \right)$$
\end{lemma}

It follows that for $q$ fixed, we get
$$\lim_{p \to \infty} p^2\; \length(\gamma)  = \frac {2\pi} {\length(m)^2}$$
where $\length(m)$ is the length of the meridian in the Euclidean cusp, normalized to have area $1$.

Putting this together with Theorem A and observing that the radius $T$ of an embedded
tube around $\gamma$ goes to infinity as $p \to \infty$ gives the estimate
$$\frac {\sqrt{2\pi}} {\length(m)} \le -3.993\pi\chi_\Q$$

Note that this estimate is most informative when $\length(m)\ll 1$, and implies that 
the genus of a hyperbolic knot can be bounded from below by a term
depending on the geometry of the cusp,
when some choice of meridian is much shorter than the longitude.

With more work, this estimate can be made quadratic in $m$. When $m$ is small, this estimate
is substantially better than the crude estimate above.

\begin{surgery_thm2}
Let $M$ be a $1$-cusped hyperbolic $3$-manifold. Normalize the Euclidean structure
on the cusp to have area $1$, and let $m$ be the shortest curve on
the cusp which is homologically essential in $M$. Then
$$-\chi_\Q \ge \frac 1 {2\pi\;\length(m)^2}$$
\end{surgery_thm2}
\begin{proof}
Let $\partial S$ denote the maximal horotorus cusp. By Jorgensen's inequality
(see \cite{Maskit}), every essential
slope on $\partial S$ has length at least $1$. It follows that if $\length(m) < 1$
on the normalized cusp, then the area of the {\em maximal} cusp
$\partial S$ is at least $1/\length(m)^2$.
If we do $p/q$ surgery with very large $p$, the area of $\partial S$ is almost unchanged.
On the other hand, an embedded solid torus of radius $T$ with core geodesic
$\gamma$  has area $2\pi \cdot \length(\gamma)\sinh(T)\cosh(T)$ and therefore
$$\area(\partial S) = \lim_{p \to \infty} 2\pi \; \length(\gamma) \sinh(T)\cosh(T) = \lim_{p \to \infty}
\frac {\pi} 2 \; \length(\gamma) e^{2T}$$
and therefore 
$$e^T \ge \frac {\sqrt{2}} {\length(m)\sqrt{\pi}} \; \length(\gamma)^{-1/2}$$
up to a multiplicative constant which goes to $1$ as $T\to +\infty$.

Using the approximations $\sinh(T) \sim e^T/2$ and $(T+1)/T \sim 1$, valid
up to a multiplicative constant which goes to $1$ as $T \to +\infty$,
it follows from Lemma~\ref{tube_estimate} that
there are quasimorphisms $\phi$ with $D(\phi) \le 2\pi$
and asymptotically satisfying 
$$\phi(a) \ge \frac {1} {\length(m)\sqrt{2\pi}}\; \length(\gamma)^{1/2}$$
Since $\scl(a) \le -\chi_\Q/2p$
there is an estimate
$$\frac {-\chi_\Q} {2p} \ge \scl(a) \ge \frac {\phi(a)} {2D(\phi)} 
\ge \length(\gamma)^{1/2} \frac 1 {2\;\length(m)(2\pi)^{3/2}}$$
Rearranging gives
$$4p^2 \; \length(\gamma) \le (2\;\length(m)(2\pi)^{3/2}\chi_\Q)^2$$
Since $\lim_{p \to \infty} p^2 \; \length(\gamma) = 2\pi\; \length(m)^{-2}$,
$$8\pi\; \length(m)^{-2} \le (2\;\length(m)(2\pi)^{3/2}\chi_\Q)^2$$
or, in other words,
$$-\chi_\Q \ge \frac {1} {2\pi \; \length(m)^2}$$
as claimed.
\end{proof}

This quadratic bound (in $\length(m)$) 
is a significant improvement over the linear bound.

\begin{problem}
Interpret the condition $\length(l) \ll 1$ topologically.
\end{problem}

\begin{remark}
The bound in Theorem B is not sharp.
It is only sensitive to the total area of the surface $S$
whereas we need only measure the area of that part of the surface contained in an
embedded tube, and moreover, only that part which is almost parallel to the core
geodesic.

A circle-packing theorem due to Boroczky \cite{Boroczky} says 
that at most $\frac 3 \pi$ of the area of
a pleated surface representing $S$ is inside the embedded tube, so one can replace the constant
$2\pi$ by $6$ above with no more work. Precisely, one gets the estimate
$$-\chi_\Q \ge \frac 1 {6 \; \length(m)^2}$$

Nathan Dunfield pointed out to me that the same estimate can be derived
from Theorem 5.1 in \cite{Agol} (moreover, the observation that Boroczky's
theorem is relevant is also made in \cite{Agol}).
\end{remark}

\begin{remark}
The most interesting thing about the geometric estimates obtained in this
section is that they give the correct
order of magnitude relationships between $\length(\gamma)$, $p$, $\length(m)$ and $\chi_\Q$.
It follows that the straightforward estimate of stable commutator length
obtained in Lemma~\ref{Seifert_estimate_scl} is typically close to
optimal.
\end{remark}

\section{Length inequalities and the Spectral Gap Theorem}\label{spectral_gap_section}

\subsection{Pleated surfaces}

Stable commutator length in the fundamental group of a manifold can be
studied by probing the manifold topologically by maps of surfaces. Under suitable geometric
hypotheses it makes sense to take representative maps of surfaces which are
tailored to the geometry of the target manifold. For $M$ hyperbolic,
the class of maps of surfaces into $M$ we use are so-called {\em pleated surfaces},
as introduced by Thurston \cite{Thurston_notes}, Chapter 8.

\begin{definition}
Let $M$ be a hyperbolic manifold. A {\em pleated surface} is a complete hyperbolic
surface $S$ of finite area, together with an isometric map $f:S \to M$ which takes
cusps to cusps, and such that every point in $S$ is in the interior of a straight
line segment in $S$ which is mapped by $f$ to a straight line segment in $M$.
\end{definition}

Note that ``isometric map'' here means that $f$ takes rectifiable curves in $S$
to rectifiable curves in $M$ of the same length. The set of points $L \subset S$
where the line segment through $p$ is unique is called the {\em pleating locus}.
The set $L$ has the structure of a {\em geodesic lamination} of $S$; i.e. a closed
union of disjoint simple geodesics. The restriction of $f$ to each component of
$S-L$ is totally geodesic. One may add geodesics to the pleating locus
if necessary so that all complementary regions are ideal triangles.

A map $f:S \to M$ has a {\em pleated representative} if there is a pleated surface
in the relative homotopy class of $f$.
If $f:S \to M$ is a pleated surface, $\area(S) = -2\pi\chi(S)$ by the
Gauss-Bonnet theorem.

\begin{lemma}\label{pleated_surface_exists}
Let $M$ be a hyperbolic manifold.
Every homotopy class of map $f:S \to M$ sending $\partial S$ to a nontrivial
hyperbolic element of $\pi_1(M)$ can either be replaced by a map $f:S' \to M$ sending
$\partial S'$ to the same element, where $\genus(S') < \genus(S)$, or it admits
a pleated representative.
\end{lemma}
For a proof, see \cite{Thurston_notes}, especially \S~8.8 and \S~8.10 or
\cite{Calegari_scl_monograph} Chapter 3.

\subsection{Length inequality}

In general, the relation between hyperbolic length and stable commutator length
is not straightforward, as the following example shows.

\begin{example}\label{example_in_surface}
Stable commutator length is characteristic, and therefore constant on orbits
of $\Aut(G)$. Let $G=\pi_1(S)$ where $S$ is a hyperbolic surface. Except for
small examples, the orbit of every conjugacy class under $\Out(G)$ is infinite.
Let $a \in \pi_1(S)$ be a conjugacy class represented by a geodesic $\gamma$ in
$S$. If $\psi \in \Aut(G)$ has a pseudo-Anosov representative, the
geodesic representatives of the conjugacy classes of elements $\psi^n(a)$
have unbounded length in $S$. On the other hand, $\scl(a) = \scl(\psi^n(a))$
is constant.
\end{example}

However the first result of this section says that an element
with {\em sufficiently small} stable commutator length necessarily
has a geodesic representative with a small length.

\begin{length_gap_thm}[Length inequality]
For every dimension $n$ and any $\epsilon > 0$ there is a constant $\delta(\epsilon,n)>0$ 
such that if $M$ is a complete hyperbolic $n$-manifold, and $a \in \pi_1(M)$
has stable commutator length $\le \delta(\epsilon,n)$ then if $\gamma$ is the
unique geodesic in the free homotopy class associated to the conjugacy class of $a$,
$$\length(\gamma) \le \epsilon$$
\end{length_gap_thm}
\begin{proof}
Let $S$ be a surface of genus $g$ with one boundary component, and
$f:S \to M$ a map wrapping $\partial S$ homotopically $m$ times around $\gamma$.
We have $\chi(S) = 1-2g$.
By Lemma~\ref{pleated_surface_exists}, after possibly reducing the genus of $S$ if
necessary, we can assume without loss of generality that $f,S$ is a pleated surface.
Note that $\area(S) = -2\pi\chi(S)$, by Gauss-Bonnet, and can be decomposed into
$-2\chi(S)$ ideal triangles. Since the boundary of $S$ wraps $m$ times around $\gamma$,
there is an equality $\length(\partial S) = m \; \length(\gamma)$.

If $\partial S$ is very long compared to the area of $S$, the surface
$S$ must become very thin along most of its boundary. Choose $\epsilon$
which is smaller than the $2$-dimensional Margulis constant $\epsilon(2)$. We defer
the precise choice of $\epsilon$ for the moment.
Consider the thick-thin
decomposition of $S$ with respect to $\epsilon$.
Note that since $S$ is a manifold with boundary, it makes sense to take the thick-thin
decomposition of the double $DS$ and then restrict this decomposition to $S$.
Each component of the thick region is relatively bounded in $S$ by proper arcs of length $\epsilon$.
Think of the thick region as a hyperbolic surface, possibly with cusps which have been ``neutered".
We distinguish between the components of $\partial C \cap \partial S$, which we call the
{\em sides} of $C$, and the components of $\partial C \cap \text{interior}(S)$ (each
of which has length $\epsilon$) which we call the {\em cusps}.
It follows that for each such component $C$ there is an inequality
$$\area(C) \ge \frac \epsilon 2 \; \length(\partial C \cap \partial S)$$
This gives a lower bound on the length of $\partial S \cap S_{<\epsilon}$.
We would like to conclude that $\partial S \cap S_{<\epsilon}$ contains a long
connected segment; to do this, we need to bound the number of components of
$\partial S \cap S_{\ge \epsilon}$.

The area of $C$ can be bounded from below from the genus of $C$, and the number of
boundary components. The worst case is when each $C$ is a neutered polygon with $n(C)$ sides, 
in which case the Euler characteristic of the corresponding un-neutered polygon is
$\chi(C) = - \frac {n(C)-2} 2$. The area of $C$ can be estimated from $\chi(C)$ by
subtracting area approximately $\epsilon$ for each neutered cusp. That is,
$$\area(C) \ge (n(C) - 2)\cdot (\pi - \epsilon)$$
From this we see that the number of
components of $\partial S \cap S_{\ge \epsilon}$ is maximized when each $C$ is a neutered
triangle, in which case there are exactly $-6\chi(S) = 12g-6$ components, and
the sum of the lengths of these components is at most $(2g -1)\pi/\epsilon$.

Therefore, by the pigeonhole principle, it follows that there is a component $\sigma$ of
$\partial S \cap S_{<\epsilon}$ of length at least
$$\length(\sigma) \ge \frac {m \cdot \length(\gamma) - (2g-1)\pi/\epsilon} {12g - 6} =
\frac {m \cdot \length(\gamma)} {12g - 6} - \frac \pi {6\epsilon}$$

By construction, $S$ contains a rectangular strip $R$ of thickness $\le \epsilon$
with $\sigma$ on one side. We denote the side opposite to $\sigma$ by $\sigma'$.
We call $\sigma$ and $\sigma'$ the {\em long} sides
of $R$. Because $S$ is oriented, the orientations on opposite
sides of $R$ are ``anti-aligned". We lift $R$ to the universal cover $\H^n$, and by abuse of 
notation refer to the lifted rectangle as $R$. The sides $\sigma,\sigma'$
of $R$ are contained in geodesics $l,l'$ which are axes for elements
$b,b' \in \pi_1(M)$ which are both in the conjugacy class $a$. Note that
the translation lengths of $b,b'$ on $l,l'$ respectively are both equal to
$\length(\gamma)$, and moreover these elements act by translating in nearly opposite directions.
We let $c \in \pi_1(M)$ 
be such that $b' = cbc^{-1}$. 

Let $p$ be the midpoint of $\sigma$, and let
$q$ be a point on the opposite side of $R$ with $d(p,q) < \epsilon$. Suppose further
that $$\length(\sigma) = \length(\sigma') > 2\;\length(\gamma)+4\epsilon$$
It follows that $b'(q) \in \sigma'$ and there is $r \in \sigma$ with
$d(b'(q),r)\le\epsilon$ and therefore $d(b'(p),r)\le 2\epsilon$. Since $qb'(q)$ and $pr$ are geodesics,
$$|d(p,r)-\length(\gamma)| \le 2\epsilon$$ and therefore $d(p,b(r)) \le 2\epsilon$ and
we can estimate
$$d(p,bb'(p)) \le 4 \epsilon$$
Similarly we estimate $d(p,b'b(p)) \le 4 \epsilon$. See Figure \ref{comparison_quad}. In
the figure, the axes $l$ and $l'$ are both roughly vertical. The element
$b$ translates points roughly downwards, and $b'$ translates points roughly upwards.

\begin{figure}[ht]
\labellist
\small\hair 2pt
\pinlabel $p$ at 110 210
\pinlabel $q$ at 45 210
\pinlabel $b'(q)$ at 20 380
\pinlabel $b'(p)$ at -30 355
\pinlabel $r$ at 120 355
\pinlabel $b(r)$ at 120 175
\pinlabel $bb'(p)$ at 220 140
\endlabellist
\centering
\includegraphics[scale=0.4]{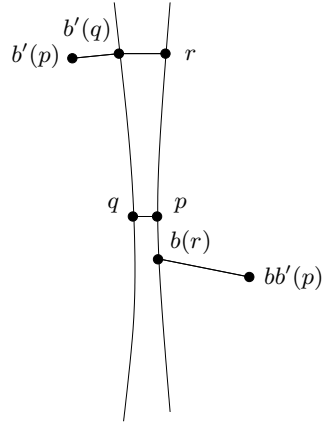}
\caption{The composition $bb'$ translates the midpoint $p$ a small distance}\label{comparison_quad}
\end{figure}

If we choose $4\epsilon$ less than an $n$-dimensional Margulis constant $\epsilon(n)$ then
$bb' = bcbc^{-1}$ and $b'b = cbc^{-1}b$ must commute. There are two possibilities, which
break up into subcases.

\vskip 12pt

{\noindent \bf Case i:} $bb'$ and $b'b$ are hyperbolic with the same
axis. In this case, since they are
conjugate, they are either equal or inverse. 

{\noindent \bf Case ia:} $bb' = b'b$. In this case $b$ and $b'$ commute, and since they
are conjugate, they are equal or inverse. But $b$ and $b'$ translate their axes in almost
opposite directions, so they cannot be equal; hence we must have $b' = b^{-1}$, and
therefore $cbc^{-1} = b^{-1}$ which implies that $c$ has order $2$, which is impossible
in a hyperbolic manifold group. 

{\noindent \bf Case ib:} $bb' = b^{-1}(b')^{-1}$. In this case $b^2 = (b')^{-2}$ and
therefore $b = (b')^{-1}$ which we have already seen is impossible in a hyperbolic manifold group.

\vskip 12pt

{\noindent \bf Case ii:} $bb'$ and $b'b$ are parabolic with the same fixed point
$z \in S^{n-1}_\infty$. In this case, $b^{-1}(bb')b$ is parabolic with fixed point
$b^{-1}(z)$. But $b^{-1}(bb')b = b'b$ which has fixed point $z$, so $b^{-1}(z) = z$.
Since $b$ is hyperbolic, we have constructed a hyperbolic and a parabolic element in
$\pi_1(M)$ with a common fixed point; this is well-known to violate discreteness, see
e.g. page 19, \cite{Maskit} for details.

\vskip 12pt

In every case we obtain a contradiction, and therefore we must have
$$2\;\length(\gamma)+4\epsilon \ge \length(\sigma)$$
Putting this together with our earlier inequality, we obtain
$$2\;\length(\gamma)+4\epsilon \ge \frac {m \; \length(\gamma)} {12g - 6} - \frac \pi {6\epsilon}$$
Rearranging this gives
$$\length(\gamma)\; \left( \frac m {12g-6} - 2 \right) \le 4\epsilon + \frac \pi {6\epsilon}$$
The right hand side is a constant which depends only on the size of
a Margulis constant in dimension $n$. If $\scl$ is sufficiently small, we can choose pleated
surfaces $S$ for which the ratio $m/(12g-6)$ is sufficiently large,
and therefore we obtain an upper bound on $\length(\gamma)$ which goes to $0$ as $\scl \to 0$
as claimed.
\end{proof}

\begin{remark}\label{optimal_estimate_remark}
By Lemma~\ref{Seifert_estimate_scl} and Theorem A, we know
that in dimension $3$ the dependence of the optimal $\delta$ on $\epsilon$ is given by
$$\delta(\epsilon,3) = O(\epsilon^{1/2})$$
for small $\epsilon$.
Similarly, in any fixed dimension $n$ there are inequalities
$$\delta(\epsilon,3) \ge \delta(\epsilon,n) \ge O(\epsilon^{(n-1)/(n+1)})$$
again for sufficiently small $\epsilon$,
but it is not clear what the optimal $\delta$ is. The
first inequality arises here since a hyperbolic $3$-manifold group acts isometrically
on $\H^n$ by stabilizing a totally geodesic $3$-dimensional subspace. The
second inequality follows from Remark~\ref{Reznikov_remark}.
\end{remark}

\begin{remark}
On the other hand, the dependence of $\delta$ on $\epsilon$ is not proper,
even in dimension $2$, as illustrated in Example~\ref{example_in_surface}.
In particular, the constant $\delta(\epsilon,n)$ is bounded above
by some universal finite bound, independent of the dimension $n$. This universal
bound is the ``homological'' Margulis constant; in the next subsection we
give an explicit estimate for this constant, and show that it is contained
between $1/12$ and $1/2$.
\end{remark}

\subsection{Spectral gap theorem}

\begin{first_limit_thm}[Spectral Gap]
Let $M$ be a closed hyperbolic manifold, of any dimension. 
Let $\delta_\infty$ be
the first accumulation point for stable commutator length on conjugacy classes. Then 
$$\frac 1 {12} \le \delta_\infty \le \frac 1 2$$
\end{first_limit_thm}
\begin{proof}
We use the same setup and notation as in the proof of Theorem~C.
Since $M$ is a closed hyperbolic manifold, there are only finitely many conjugacy classes
represented by geodesics shorter than any given length. So we suppose $a$
is a conjugacy class represented by a geodesic $\gamma$ which is ``sufficiently
long" (in a sense to be made precise in a moment). We choose $\epsilon$
and find a segment $\sigma$, as in the proof of Theorem~C,
and suppose we have
$$\length(\gamma) + 4\epsilon < \length(\sigma)$$
(note the missing factor of $2$).
We choose $p$ to be one of the endpoints of $\sigma$, so that
$$d(p,bb'(p)) \le 4 \epsilon$$

Since $M$ is fixed, there is some $\epsilon$ such that $4\epsilon$ is smaller than
the translation length of any nontrivial element in $\pi_1(M)$. Hence $bb' = \id$.
But this means $cbc^{-1} = b^{-1}$, and $c$ has order $2$, which is impossible
in a hyperbolic manifold group.

Contrapositively, this means that we must have
$$\length(\gamma) + 4\epsilon \ge \length(\sigma)$$
and therefore, just as in the proof of Theorem~C,
we obtain
$$\length(\gamma)\cdot \left( \frac m {12g-6} - 1 \right) \le 4\epsilon + \frac \pi {6\epsilon}$$
In contrast to the case of Theorem~C, the
right hand side definitely depends on the manifold $M$. Nevertheless, for fixed $M$,
it is a constant, and we see that for $\gamma$ sufficiently long, $g/m$ cannot be
much smaller than $1/12$. This establishes the lower bound in the theorem.

\vskip 12pt

To establish the upper bound, observe that since $M$ is a closed hyperbolic
manifold, $\pi_1(M)$ contains many nonabelian free groups. In fact, if $a,b$ are
arbitrary noncommuting elements of $\pi_1(M)$, sufficiently high powers of
$a$ and $b$ generate a free group, by the ping-pong lemma.
This copy of $F_2$ is quasi-isometrically embedded, and by passing to a subgroup,
one obtains quasi-isometrically embedded copies of free groups of any rank.

Let $S$ be a once-punctured hyperbolic surface of genus $2$. Then $\pi_1(S)$ is free
of rank $4$. Let $\phi:\pi_1(S) \to \pi_1(M)$ be a quasi-isometric injection,
as above. Let $T$ be a once-punctured torus, embedded as a subsurface of $S$,
and let $a$ be the conjugacy class in $\pi_1(S)$ corresponding to $\partial T$.
By construction, $\scl(a) \le 1/2$. Moreover, under
the action of the mapping class group of $S$, the orbit of $a$ is
infinite. So there are infinitely many conjugacy classes $a_i$ in $\pi_1(S)$
with $\scl(a_i) \le 1/2$. If we represent the $a_i$ by geodesics in $S$,
their lengths grow without bound. Since $\phi$ is a quasi-isometric embedding,
the geodesics in $M$ corresponding to the elements $\phi(a_i)$ have unbounded
length, and therefore they correspond to infinitely many conjugacy
classes in $\pi_1(M)$. Since $\scl$ is monotone under homomorphisms, we have
constructed infinitely many conjugacy classes in $\pi_1(M)$ with $\scl \le 1/2$.
\end{proof}

\begin{remark}
If $M$ is a finite volume hyperbolic $3$-manifold, possibly with cusps,
the first accumulation point $\delta_\infty$ still satisfies $\delta_\infty \ge 1/12$.
For, if $\gamma$ is a sufficiently long closed geodesic in $M$, we can find a Dehn filling $M'$
of $M$ in which the image of $\gamma$ is also homotopic to an arbitrarily long
closed geodesic. The stable commutator length of a sufficiently long geodesic 
in $M'$ is at least $1/12-\epsilon$ for any positive $\epsilon$, by Theorem D. Moreover, 
the stable commutator length can only go down under Dehn filling, since
stable commutator length is non-increasing under homomorphisms; the claim follows.

If $M$ is a finite volume hyperbolic $n$-manifold, the argument of Theorem C shows
directly that $\delta_\infty \ge 1/24$. But it seems reasonable to 
expect that $\delta_\infty \ge 1/12$ in this case too.
\end{remark}

\begin{remark}
It is a well-known conjecture that there is a uniform positive
lower bound on the length of a closed geodesic in an arithmetic hyperbolic
manifold. This is related to Lehmer's conjecture in number theory; 
c.f. \cite{Lehmer},\cite{Neumann_Reid}.
This motivates the following

\begin{conjecture}
For each $n$ there is a constant $C(n)>0$ such that for any arithmetic hyperbolic $n$-manifold $M$
and any non-parabolic $a \in \pi_1(M)\backslash \id$, there is an inequality
$$\scl(a) \ge C(n)$$
\end{conjecture}

Compare with the following Example, due to Agol.
\begin{example}[Agol]
In \cite{Agol_systole}, Agol shows that in any dimension $n$, if there is an arithmetic
hyperbolic $n$-manifold defined by a quadratic form which is GFERF (i.e. subgroup separable
on geometrically finite subgroups) then there exist finite volume hyperbolic $n$-manifolds
with arbitrarily short geodesics. Such subgroups exist for all $n \le 8$.
The examples are obtained from arithmetic manifolds with lots of
immersed totally geodesic hypersurfaces by mutation.
\end{example}
\end{remark}

\begin{remark}\label{inverse_conjugate_vanishes}
If ``manifold'' is replaced by ``orbifold'' in the statement of Theorem C or Theorem D, 
then there are a number of places where the argument must be modified. In an orbifold,
elliptic elements with a common fixed point do not necessarily commute. Moreover,
we must allow the possibility that $c^2=\id$ in which case $b$ is conjugate
to $b^{-1}$.

In any group, if there is $c$ such that $cb^{-1}c^{-1} = b$ then 
$b^{2n} = b^ncb^{-n}c^{-1} = [b^n,c]$
and the stable commutator length of $b$ vanishes. In particular,
an element $a$ in a hyperbolic orbifold group is either conjugate to its
inverse or satisfies $\scl(a)>0$ where the bound may be estimated in terms
of the translation length of $a$ and the order of torsion elements in the group.
\end{remark}

\begin{remark}
Theorem C and Theorem D admit natural generalizations to (word-) 
hyperbolic groups, and to groups acting on $\delta$-hyperbolic spaces, 
such as mapping class groups and amalgamated free products. 
These results will appear in a forthcoming paper \cite{Calegari_Fujiwara}.
\end{remark}

\subsection{Broken Windows}

The following discussion is included for the benefit of the reader, and is not logically
necessary for the rest of the paper.

If an $n$th power of an element in $\pi_1(M)$
can be expressed as a product of $m$ commutators,
one obtains an immersed $2$-complex $K \subset M$ obtained from wrapping the
boundary of a genus $m$ surface $n$ times around the corresponding geodesic
$\gamma$. The
map of $K$ into $M$ is almost never $\pi_1$-injective, but when it is, one
can derive a global relationship between $n,m$ and $\length(\gamma)$.

\begin{theorem}[Thurston]\label{injective_wrap_estimate}
Let $n \ge 3$ and let $G_{n,m}$ be the group defined in terms of generators and relations by
$$G_{n,m} = \langle a,b_1,c_1,\dots,b_m,c_m \; | \; a^n = \prod_i [b_i,c_i] \rangle$$
Let $M$ be a complete hyperbolic $3$-manifold with $\pi_1(M) = G_{n,m}$. Then
there is a constant $C>0$ such that if $\gamma$ is the geodesic in $M$ corresponding to
the conjugacy class $a \in G_{n,m}$, there is an estimate
$$\length(\gamma) \le C \frac m {\log(n-1)}$$
\end{theorem}
This is actually an effective version of a special case of Thurston's 
``Broken Windows Only Theorem". We give a sketch of a proof; for details,
consult \cite{Thurston_III}.
\begin{proof}
Let $K \subset M$ be an embedded $2$-complex built out of the image of a suitable
pleated surface $S$ by wrapping the boundary $m$ times around the geodesic
$\gamma$. By Gauss-Bonnet,
$$\area(K) = (4m-2)\pi$$
In the universal cover $\til{M} = \H^3$ the preimage $\til{K}$ is a tree of
hyperbolic planes in its intrinsic metric, with valence $m\ge 3$ along the
branch locus, which is equal to the preimage $\til{\gamma}$ of $\gamma$. A random walk
on $\til{K}$ pushes down to a random walk on $K$. Since the geodesic curvature
of $K$ is uniformly negative, a random walk crosses $\gamma$ with frequency
proportional to the length of $\gamma$.

It follows that there is a constant $C>0$ independent of $m$ and $n$
such that for a random geodesic $\sigma$ in $\til{K}$ of length $T$, the
number of intersections of $\sigma$ with $\til{\gamma}$ satisfies
$$\# \lbrace \sigma \cap \til{\gamma} \rbrace = \frac {CT\cdot \length(\gamma)} m + o(T)$$
Since $\til{K}$ has valence $n$ along each such sheet, this implies that
the area of the ball of radius $T$ about a point $p \in \til{K}$ satisfies
$$\area(B_{\til{K}}(T,p)) \gtrsim (n-1)^{CT\cdot \length(\gamma)/m}\cdot e^T$$

On the other hand, since $K$ is immersed in $M$, this growth rate must be
less than the growth rate of the volume of hyperbolic $3$-space, 
or else $\til{K}$ would accumulate on itself. It follows that we can estimate
$$T + \frac {\log(n-1)CT\cdot\length(\gamma)} m \le 2T + \text{constant}$$
and we are done.
\end{proof}

\begin{example}
Let $S$ be the connect sum of three projective planes. A presentation for
$\pi_1(S)$ is
$$\pi_1(S) = \langle a,b,c \; | \; [a,b] = c^2 \rangle$$
The element $c$ is represented by a simple one-sided
geodesic $\gamma$ on $S$. Geometrically, $S$ is obtained from a once-punctured
torus with geodesic boundary by mapping the boundary to $\gamma$ by a double
covering. A hyperbolic once-punctured torus exists with boundary of any given
positive length; it follows that the length of $\gamma$ may be specified arbitrarily.
This example shows that the restriction $n\ge 3$ is necessary
in Theorem~\ref{injective_wrap_estimate}.
\end{example}

\section{de Rham length versus stable commutator length}\label{de_Rham_section}

\subsection{de Rham length}

The results in \S~\ref{surgery_section} show that de Rham quasimorphisms give good
lower bounds on $\scl$ for short geodesics contained in embedded tubes.
It is natural to ask how rare de Rham quasimorphisms are amongst all
quasimorphisms in general.

\begin{definition}
Let $M$ be a hyperbolic manifold, and let $a \in \pi_1(M)$. Let $\gamma$ be
the unique geodesic in $M$ in the free homotopy class of the conjugacy class
of $a$. The {\em de Rham length} of $a$, denoted $\drl(a)$, is the supremum
$$\drl(a) = \sup_{\alpha \in \Omega^1} \frac {\int_\gamma \alpha} {2\pi\|d\alpha\|}$$
\end{definition}

By Theorem~\ref{Bavard_theorem} and Lemma~\ref{form_is_quasimorphism}
there is an inequality $\scl(a) \ge \drl(a)$ for any $a\in\pi_1(M)$.

\begin{remark}
de Rham length is naturally related to norms from geometric measure
theory, especially the {\em Whitney flat norm} (see \cite{Whitney},\cite{Federer}).
\end{remark}

\subsection{Ideal Quadrilaterals}

The following construction is due to Thurston.

\begin{example}[Thurston]\label{small_area_example}
Let $P$ be an ideal quadrilateral in $\H^n$, $n>2$ with the property that opposite vertices
are very close together in $S^{n-1}_\infty$, and the ideal simplex they span has
simplex parameter whose imaginary part is much bigger than its real part (this requires $n>2$).
Then $P$ spans a disk with arbitrarily small area, made from a bisecting ``square''
$S$ together with four ``cusps'' $C_1,\dots,C_4$. See Figure~\ref{small_disk}.

\begin{figure}[ht]
\labellist
\small\hair 2pt
\pinlabel $C_1$ at 103 210
\pinlabel $C_3$ at 100 70
\pinlabel $C_4$ at 215 105
\pinlabel $C_2$ at 217 190
\pinlabel $S$ at 108 140 
\endlabellist
\centering
\includegraphics[scale=0.5]{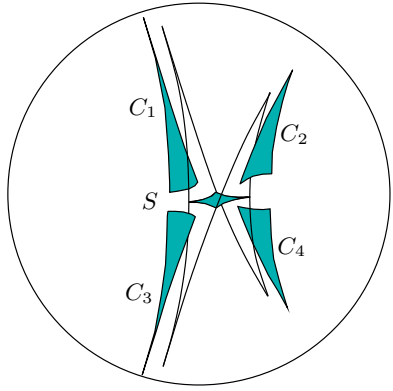}
\caption{A disk made from $S$ and the $C_i$ spans $P$}\label{small_disk}
\end{figure}

For any finite volume hyperbolic $n$-manifold $M$ we can immerse $P$ in $M$ by
projecting $\H^n \to M$. Now, for generic choices of endpoints, the four boundary
geodesics are almost perfectly equidistributed in $UTM$, and therefore after an arbitrarily
small perturbation, we can glue up opposite faces to product a once-punctured
torus $T$ with arbitrarily small area. For generic translation lengths
along the sides, the hyperbolic structure on the torus is incomplete, and
can be completed by adding a geodesic boundary $\gamma$.

Since the matching can be performed with infinitely many
translation lengths along the sides, $\length(\gamma)$ can be chosen to
be arbitrarily long.
\end{example}

This construction produces infinitely many closed geodesics $\gamma$ in $M$
bounding surfaces of arbitrarily small area. By Stokes' theorem, if
$a \in \pi_1(M)$ is a conjugacy class corresponding to such a $\gamma$, then
$\drl(a)$ is arbitrarily small.

On the other hand, by Theorem C, the stable commutator length $\scl(a)$ 
cannot be too small, when $\length(\gamma)$ is large. This proves the following
theorem.

\begin{bad_comparison_thm}
Let $M$ be a finite volume complete hyperbolic $n$-manifold with $n>2$. Then for
any $\epsilon >0$ there are infinitely many conjugacy classes $a \in \pi_1(M)$ with
$$\drl(a) \le \epsilon\; \scl(a)$$
\end{bad_comparison_thm}

\begin{problem}
Give an algebraic characterization of conjugacy classes $a \in \pi_1(M)$
with small $\drl$.
\end{problem}

\section{Stable commutator length in other $3$-dimensional geometries}\label{geometry_section}

For completeness, we discuss lower bounds on $\scl$ in 
$\pi_1(M)$ for $M$ a $3$-manifold with one
of the seven (non-hyperbolic) $3$-dimensional geometries.

\begin{example}[$S^3,S^2\times \R,\R^3, \Sol,\Nil$]
If $M$ has one of these geometries, $\pi_1(M)$ is amenable, and therefore
by Example~\ref{amenable_example_vanishes}, $\scl$ vanishes identically
on $[\pi_1(M),\pi_1(M)]$.

In most cases, one can see $\scl(a)=0$ directly for
trivial reasons, but the case that $M$ has a $\Sol$ geometry is a little more interesting,
and merits discussion.

A typical Sol manifold $M$ has a fundamental
group $\pi_1(M) = \Z^2 \rtimes \Z$ where the conjugation action
of the generator $g$ of $\Z$ on $\Z^2$ is given by an Anosov matrix $A$
(i.e. one with $|\text{trace}|>2$). In this case, the kernel of the map
$\pi_1(M) \to H_1(M;\R)$ is equal to the normal $\Z^2$ subgroup.

So suppose $a \in \Z^2$. Let $v_1,v_2$ be two independent 
eigenvectors of $A$ over $\R^2$. 
We can express $a$ as a product of elements in $\Z^2$
$$a = w_1w_2v$$
where $w_1$ is almost a multiple of $v_1$, $w_2$ is almost a multiple of
$v_2$, and there is an inequality $|v|/|a|\le C<1$ for some constant $C$. 
For suitable non-negative integers $k_1,k_2$,
the conjugates $g^{k_1}w_1g^{-k_1}$ and $g^{-k_2}w_2g^{k_2}$ 
are both words of uniformly bounded length in $\Z^2$ with respect to any fixed generating
set, so $a$ can be written as a product of a bounded number of commutators
times something of size definitely smaller than that of $a$. 
By induction, we obtain an estimate $\cl(a) \le O(\log(|a|))$, and therefore
$\scl(a)=0$.
\end{example}

\begin{example}[$\H^2 \times \R$]
The typical case is that $\pi_1(M) = G \times \Z$ where $G$
is a hyperbolic surface group. 
If $a \in \pi_1(M)$ has nontrivial image in the $\Z$
factor, then $\scl(a) = \infty$. Otherwise, $a$ is contained in $G$.
In an orientable hyperbolic surface group, $\scl(a) \ge 1/2$ by Example~\ref{free_example_bound}.
\end{example}

\begin{example}[$\til{\SL(2,\R)}$]
The typical case is that there is a short exact sequence 
$$\Z \to \pi_1(M) \to G$$ 
where $G$ is a hyperbolic surface group. In an orientable hyperbolic surface group,
$\scl\ge 1/2$ by Example~\ref{free_example_bound}.

There is also an action of $\pi_1(M)$ on $\R$ which comes from the universal
covering $\R \to S^1$ lifting the natural action of $\PSL(2,\R)$ on $\RP^1 = S^1$.
For this action, Poincar\'e's {\em rotation number} (Example~\ref{rotation_example})
is a natural homogeneous quasimorphism with defect $1$. 
This quasimorphism takes the value $1$ on the generator
of the center of $\pi_1(M)$, and takes some value $p/q$ where $p \ne 0$ is coprime to
$q$ on any element which maps to a torsion element of order $q$ in $G$.
\end{example}

\end{document}